\documentclass{amsart}
\usepackage{graphicx}
\usepackage{latexsym}
\usepackage{amsfonts}
\usepackage[all]{xy}
\usepackage{amssymb, mathrsfs, amsfonts, amsmath}
\usepackage{amsbsy}
\usepackage{amsfonts}
\usepackage{amssymb}
\newtheorem{theorem}{Theorem}[section]

\theoremstyle{definition}

\theoremstyle{remark}
\newtheorem{remark}[theorem]{Remark}

\numberwithin{equation}{section}
\newcommand{\hess}{\textrm{Hess}\,}
\newcommand{\grad}{\textrm{grad}\,}

\newcommand{\Hvec}{\Vert \hspace{-1mm}\stackrel{\to}{H}\hspace{-1mm}\Vert}

\begin{document}

\title{On the mean curvature\\ of  Nash isometric embeddings}
\author{G. Pacelli Bessa \and J. Fabio Montenegro  }
\date{\today}
\maketitle

\begin{abstract}
{\noindent J. Nash proved in  \cite{nash} that the geometry of any Riemannian manifold $M$ imposes no restrictions to be embedded isometrically into a (fixed)  ball $B_{\mathbb{R}^{N}}(1)$ of  the Euclidean space $\mathbb{R}^{N}$. However,  the geometry  of $M$ appears, to some extent, imposing  restrictions on the mean curvature vector of the embedding.   }
\end{abstract}
\section{Introduction}
In 1956, John Nash in a celebrated article \cite{nash}, proved that any complete $n$-dimensional Riemannian manifold $M$ can be isometrically embedded in a ball $B_{\mathbb{R}^{N}}(1)$  of the Euclidean space $\mathbb{R}^{N}$  of  radius $1$ with $N=n(n+1)(3n+11)/2$. Although the  geometry of $M$ does not impose  restrictions on the existence of an isometric embedding $\varphi:M\hookrightarrow B_{\mathbb{R}^{N}}(1)$, it does impose restrictions on  the mean curvature vector of $\varphi$. That is
  expressed in the following theorem.
\begin{theorem}Let $\varphi :M \hookrightarrow B_{\mathbb{R}^{N}}(R)\subset \mathbb{R}^{N}$ be an isometric embedding of a complete $n$-dimensional Riemannian manifold into the Euclidean space  $\mathbb{R}^{N}$. Then  $$\sup_{M}\Hvec \geq n/R-\lambda^{\ast}(M)\cdot R/2$$
\label{thm1}Where $\lambda^{\ast}(M)$ is the fundamental tone of $M$  given by
$$
\lambda^{\ast}(M )=
\inf\left\{\frac{\int_{M}\vert \nabla f\vert^{2}}{\int_{M} f^{2}},\,f\in
{ H_{0}^{1}(M )\setminus\{0\}}\right\},
$$$\stackrel{\to}{H}={\rm tr}\alpha$ is the mean curvature vector and $\alpha$ is the second fundamental form of $\varphi$.
\end{theorem}
This result or rather its proof   was already  given in the literature, for instance in \cite{bessa-montenegro}, \cite{bessa-montenegro2}, \cite{bessa-silvana}, \cite{jorge-xavier}  but these   restrictions on the mean curvature vector of  Nash embeddings was never observed.
\section{Proof of Theorem \ref{thm1}} Theorem \ref{thm1} is an application of the following variation of  Barta's Theorem \cite{barta}.
\begin{theorem}Let $M$ be a complete Riemannian manifold.  Let $f$ be a positive smooth function on $M$. Then $\lambda^{\ast}(M)\geq \inf_{M}(-\triangle f/f)$.
\end{theorem} The proof of Theorem \ref{thm1} is as follows.
Let $\varphi:M\hookrightarrow B_{\mathbb{R}^{N}}(1)\subset\mathbb{R}^{N}$ be an isometric immersion. Let $g:\mathbb{R}^{N}\to \mathbb{R}$ defined by $g=(R^{2}-\rho^{2})/2$, where $\rho$ is the distance function in $\mathbb{R}^{N}$ to the origin. Define $f=g\circ \varphi:M\to \mathbb{R}$. Observe that $f>0$ thus by Barta's Theorem we have that $\lambda^{\ast}(M)\geq \inf_{M}(-\triangle f/f)$. To compute $\triangle f$ we use the following formula proved by Jorge-Koutrofiotis in \cite{jorge}.
\begin{equation}\triangle f (y)=\sum_{i=1}^{n}\hess g (\varphi (y)(e_{i},e_{i})+ \langle \grad g, \stackrel{\to}{H}\rangle (\varphi (y))\end{equation}where $\{e_{i}\}$ is an orthonormal basis of $T_{y}M$ identified with $\{d\varphi (y) e_{i}\}$.

 \vspace{2mm}
 \noindent Choose a basis  $\{e_{i}\}$ as follows. Let $\{\grad \rho, \partial/\partial \theta_{1}, \ldots \partial/\partial_{N-1}\}$ be a polar basis for $T_{\varphi (y)}\mathbb{R}^{N}$ and
  let $e_{1}=\langle e_{1}, \grad \rho\rangle \grad \rho + \langle e_{1}, \grad \rho^{\perp}\rangle \grad \rho^{\perp}$ and $e_{2},\ldots, e_{n}\in \{ \partial /\partial \theta_{1}, \ldots \partial/\partial_{N-1}\}$.
With this choice  we obtain
 \begin{eqnarray}\triangle f &=& \sum_{i=1}^{n}\hess g (e_{i},e_{i})+ \langle \grad g, \stackrel{\to}{H}\rangle \nonumber\\
 & =& g''\langle e_{1}, \grad \rho  \rangle^{2}+\frac{g'}{\rho}\langle e_{1}, \grad \rho^{\perp}  \rangle^{2}+\frac{m-1}{\rho}g'+ g'\langle \grad \rho, \stackrel{\to}{H}\rangle \nonumber \\
 &=&( g''-\frac{g'}{\rho})\langle e_{1}, \grad \rho  \rangle^{2}+ \frac{m}{\rho}g'+ g'\langle \grad \rho, \stackrel{\to}{H}\rangle \nonumber\\
 &=& -n-\rho \langle \grad \rho, \stackrel{\to}{H}\rangle .\nonumber
   \end{eqnarray}  From this we have
\begin{equation}-\frac{\triangle f}{f}= \frac{2n}{R^{2}-\rho^{2}}+\frac{2\rho}{R^{2}-\rho^{2}}\langle \grad \rho, \stackrel{\to}{H}\rangle \geq \frac{2n}{R^{2}}-  \frac{2\sup_{M}\Hvec }{R}. \end{equation} We may assume that $\sup_{M}\Hvec <\infty$, otherwise there is nothing to prove. We have then $$\lambda^{\ast}(M)\geq \inf ( -\triangle f/f)\geq  \frac{2n}{R^{2}}-  \frac{2\sup_{M}\Hvec }{R}$$Therefore \begin{equation}\sup_{M}\Hvec\geq n/R-\lambda^{\ast}(M)\cdot R/2.\label{eqThm}\end{equation}

 \noindent If we consider the mean curvature function $H$ defined by $\stackrel{\to}{H}=n\cdot H\cdot \stackrel{\to}{\eta}$, where $\stackrel{\to}{\eta}$ is unit vector normal to the submanifold, the inequality (\ref{eqThm}) becomes $\sup_{M}\vert H \vert \geq 1/R-\lambda^{\ast}(M)\cdot R/2n$.
\begin{remark}The hyperbolic $n$-space $\mathbb{H}^{n}(-1)$ has  $\lambda^{\ast}(\mathbb{H}^{n}(-1))=(n-1)^2/4$. Thus it can not be embedded (immersed)  minimally into a ball $B_{\mathbb{R}^{N(n)}}(1)$ for $n=2,\ldots, 9$.
\end{remark}
 \vspace{5mm}
\noindent{\bf Acknowledgments:}  We would like to thank Paolo Piccione for the discussions we had with respect to this note when he visited us in Fortaleza.


\begin{thebibliography}{abcd}

\bibitem{barta} J. Barta,  {\em Sur la vibration fundamental d'une membrane.} C. R. Acad. Sci. \textbf{204}, (1937), 472--473.
%

    %
\bibitem{bessa-montenegro}G. P. Bessa \and J. Fabio Montenegro,
{\em Eigenvalue estimates for submanifolds  with locally bounded
mean  curvature. } Ann.  Global Anal. and Geom.  24, 2003,
279-290.

\bibitem{bessa-montenegro2} G. P. Bessa, \and J. Fabio Montenegro, {\em An Extension of Barta's Theorem and Geometric Applications}.  Ann. Global Anal. Geom. \textbf{31}  (2007),  no. 4, 345--362.

\bibitem{bessa-silvana}G. P. Bessa, \and M. Silvana Costa, {\em     Eigenvalue Estimates for  submanifolds   with locally bounded mean curvature in  $N \times \mathbb{R}$.}

\bibitem{jorge} L. P. Jorge, \and D. Koutroufiotis, {\em An estimative for the curvature of bounded submanifolds}, Amer. J. of Math. \textbf{103}, (1980), 711--725.

\bibitem{jorge-xavier} L. P. Jorge \and F. Xavier, {\em An inequality
  between the exterior diameter and the mean curvature of bounded
  immersions.} Math. Z. \textbf{178}, (1981) 77--82.

\bibitem{nash} J. Nash, {\em The imbedding problem for Riemannian manifolds. } Ann. of Math. (2) \textbf{63},(1956) 20--63.


\end{thebibliography}
\end{document}